\documentstyle[10pt]{article}

\newcommand\qed{\nopagebreak[4]\begin{flushright}\rule{0.1in}{0.1in}
\end{flushright}\pagebreak[2]}
\newcommand\chargp[1]{\widehat{#1}}

\newcommand\Zar[1]{\overline{#1}^{\mathrm{Zar}}}

\def\Spec{{\mathrm{Spec}}}
\def\sign{{\mathrm{sign}}}

\def\lcm{\mathop{{\mathrm{lcm}}}}

\def\Epsilon{\varepsilon}
\def\OEpsilon{\overline{\Epsilon}/\Epsilon}

\def\Tor{{\mathrm{Tor}}}
\def\Hom{{\mathrm{Hom}}}

\def\rank{{\mathrm{rank}}}
\def\iso{\cong}
\def\length{{\mathrm len}}

\def\im{{\mathrm{im}}}

\def\mod{{\mathrm{mod}}}

\def\ord{{\mathrm{ord}}}
\def\per{{\mathrm{per}}}
\def\order{{\mathrm{ord}}}

\newfam\meuffam
\font\tenmeuf=eufm10
\font\sevenmeuf=eufm7
\font\fivemeuf=eufm5

\textfont\meuffam=\tenmeuf
\scriptfont\meuffam=\sevenmeuf
\scriptscriptfont\meuffam=\fivemeuf

\newfam\msbfam
\font\tenmsb=msbm10
\font\sevenmsb=msbm7
\font\fivemsb=msbm5

\textfont\msbfam=\tenmsb
\scriptfont\msbfam=\sevenmsb
\scriptscriptfont\msbfam=\fivemsb

\def\Bbb{\fam\msbfam\tenmsb}

\def\l{{\germ l}}

\def\C{{\Bbb C}}

\def\N{{\Bbb N}}

\def\Q{{\Bbb Q}}

\def\Z{{\Bbb Z}}

\def\sA{{\cal A}}

\def\sC{{\cal C}}

\def\sF{{\cal F}}

\def\sL{{\cal L}}

\def\sP{{\cal P}}

\def\sS{{\cal S}}
\def\sT{{\cal T}}
\def\sU{{\cal U}}
\def\sV{{\cal V}}

\font\l=cmr10 at 10pt
\font\ls=cmr7
\font\lss=cmr5
\font\lsy=cmsy10
\font\lsys=cmsy7
\font\lsyss=cmsy5
\font\lmi=cmmi10
\font\lmis=cmmi7
\font\lmiss=cmmi5
\font\lex=cmex10

 at 6truept
 at 6truept

\def\mapright#1{\smash{
  \mathop{\longrightarrow}\limits^{#1}}}

\newcommand\cd[1]{\matrix{#1}}

\input epsf

\newbox\figbox%
\newdimen\fight%
\def\figset#1{\setbox\figbox=\hbox{\epsffile{#1}}%
\fight=\ht\figbox\advance\fight by 1cm%
\vbox to\fight{\vfill}\box\figbox}

\newtheorem{theorem}{Theorem}
\newtheorem{lemma}{Lemma}[section]
\newtheorem{proposition}[lemma]{Proposition}

\newtheorem{corollary}[lemma]{Corollary}

\newcommand\heading[1]{\smallskip\noindent{\bf
#1}}

\title{Torsion Points on an Algebraic Subset of an Affine Torus}
\author{Eriko Hironaka\thanks{Research partially supported by
N.S.E.R.C. grant OGP0170260}}
\begin{document}
\maketitle
\begin{abstract}  
Work of Laurent and Sarnak, following a conjecture of Lang, shows
that the number of torsion points of order n on an algebraic subset 
of an affine complex torus is polynomial periodic.  In this paper, 
we find bounds on the degree and period of this number as a function 
of n.  Some examples, including the number of n torsion 
points on Fermat curves, are computed to illustrate the methods.
\end{abstract}

\section{Introduction.}
Let $V \subset (\C^*)^r$ be an algebraic subset.
Laurent has shown that the torsion points on $V$ lie on a finite
union of translates of affine subtori contained in $V$.  It follows
that the number, $p_V(n)$, of
torsion points of order $n$ on $V$ is polynomial periodic.
In this paper, we find formulas and bounds for the degree and period
of $p_V(n)$ in terms of defining equations for $V$.

Counting torsion points on algebraic subsets of an the affine torus
is useful for studying abelian representations of a finitely
presented group $\Gamma$.  For example, the Alexander invariants
define a stratification of the character variety for $\Gamma$,
which naturally embeds in $(\C^*)^r$,
and the first Betti number of unbranched coverings can be computed
from the torsion points on these strata (cf. \cite{A-S:Betti},
\cite{Hiro:Alex}).

In \cite{Laur:Equ}, Laurent gives the following description
of the set of torsion points $\Tor(V)$ on $V$.

\begin{theorem} (Laurent) 
For any algebraic subset $V$ of $(\C^*)^r$ there is a finite set
of rational planes $Q_1,..,Q_\ell \subset V$ such that 
$$
\Tor(V) =  \bigcup_{i=1}^\ell \Tor(Q_i).
$$
\end{theorem}
Here, a {\it rational plane} $Q \subset (\C^*)^r$ is a subset of the form
$\eta P$, where $P \subset (\C^*)^r$
is an {\it affine subtorus}, or connected algebraic subgroup, 
and $\eta \in (\C^*)^r$ is an element of finite order.
Laurent's result extends further than the statement we give here and
settles a more general conjecture of Lang (\cite{Lang:Conj}, p.220).

It follows from Theorem 1 (cf. \cite{A-S:Betti}) 
that $p_V(n)$ is a {\it polynomial periodic} function in $n$,
that is, there exist periodic functions $a_0(n), \dots, a_d(n)$ such
that
$$
p_V(n) = a_0(n) + a_1(n)n + \dots + a_d(n)n^d.
$$
For, by Theorem 1, the Zariski closure of $V$ has a decomposition
into a finite union of rational planes
$$
\Zar{\Tor(V)} = Q_1 \cup \dots \cup Q_\ell
$$
and hence $p_V(n)$ is given by the following formula
\begin{eqnarray}
p_V(n) = \sum_{k=1}^\ell \quad  
\sum_{1\leq i_1 < \dots < 
i_k \leq\ell}
(-1)^{\ell-k} 
p_{Q_{i_1} \cap \dots \cap Q_{i_k}}(n).
\end{eqnarray}
It is not hard to see that a finite intersection of rational
planes is  a finite union of disjoint rational planes (see Prop. 3.6
for a more precise description).   Furthermore, the number of
$n$-torsion points on a rational plane $Q$ is given by
\begin{eqnarray*}
p_Q(n) &=& \left\{\begin{array}{ll} 
n^{\dim (Q)} &\quad\mbox{if $\ord(Q)\ |\ n$,}\\
0 &\quad\mbox{otherwise,}
\end{array}\right .
\end{eqnarray*}
where $\ord(Q)$ is the least integer $n$ such that $Q^n$ contains the identity
element of
$(\C^*)^r$.

We will concentrate on finding the degree and period of $p_V$.
The {\it degree} of $p_V$, written $\deg (p_V)$, is
the largest $d$ such 
that $a_d(n)$ is not constantly
zero and 
the {\it period} of $p_V$, written $\per (p_V)$, 
is the least common multiple of the periods
of $a_0(n),\dots,a_d(n)$.  Thus, for example, if $Q \subset (\C^*)^r$ is a
rational plane, then $\ord (Q) = \per (p_Q)$.
As with ordinary polynomials the degree of $p_V$ determines the order of
growth of $p_V(n)$:
if $\deg(p_V) = d$, then
$$
p_V(n) = O(n^d)
$$
and for some fixed integer $c$
$$
p_V(n) \asymp n^d
$$
for all $n \equiv c\ (\mod\ \per(p_V))$.

If $\Zar{\Tor(V)} = Q_1 \cup \dots \cup Q_\ell$, where $Q_1,\dots,Q_\ell$
are rational planes, then (1) implies that
\begin{description}
\item{(i)} $\deg (p_V)  = \max\  \{ \dim(Q_i) : i=1,\dots,\ell \} $, 
\item{(ii)} $\per (p_V)$ divides the least common multiple of 
$\ord(Q)$ where $Q$ ranges among connected components of
$$
Q_{i_1} \cap \dots \cap Q_{i_k}, \quad 
1 \leq i_1 < \dots < i_k \leq \ell.
$$
\end{description}
In section 2, Theorem 2, we give bounds for $\deg(p_V)$
and $\per(p_V)$, when $V$ is defined over $\Q$.
More precise formulas for $\deg(p_V)$ and $\per(p_V)$, for general $V$,
are obtained later in section 5 after developing notation and 
theory in sections 3 and 4.

The main ingredients of Laurent's proof of Theorem 1 
can be stated as follows.  Given an algebraic subset
$V \subset (\C^*)^r$, one defines a finite set $\Pi$ 
of  mappings,
$$
\phi: (\C^*)^r \rightarrow \Gamma_{\phi},  \qquad \phi \in \Pi
$$
to algebraic groups $\Gamma_{\phi}$.  These mappings have the 
following properties.
\begin{description}
\item{(A)} All the fibers of $\phi$ are finite 
unions of rational planes; and
\item{(B)} 
$$
\Tor(V) = \bigcup_{\phi\in \Pi} \phi^{-1}(\sS_{\phi}),
$$
where $\sS_{\phi} \subset \Gamma_{\phi}$ are finite subsets.
\end{description}
We will go further in this paper by finding the fibers in 
(A) and the sets $S_\phi$ in (B)
explicitly and relating them to the degree and period of $p_V$.

Techniques for finding the fibers in (A) are given in section 3, 
where we develop some tools for studying rational planes and give
properties of monomial mappings.  Our approach 
is to focus on the relationship between algebraic subgroups of 
$(\C^*)^r$ and subgroups of $\Z^r$.  We use this correspondence  
to describe the rational planes in fibers of 
monomial mappings (see Lemma 3.4  and Cor. 3.5) and, as 
a consequence, in the intersection
of a finite set of rational planes (see Prop. 3.6). 

In many natural applications, for example, the Alexander strata
mentioned above, the algebraic subsets $V$ are defined over $\Q$, so we
will concentrate on this case.  
Then, the sets $S_\phi$ in (B) 
are related to formal $\Q$-linear combinations of roots of
unity.  
In section 4, we review some of the theory of $\Q$-linear relations
among roots of unity using ideas of Schoenberg.  The main idea is to
view formal $\Q$-linear relations as convex polygons with rational sides
and angles.  Schoenberg shows in \cite{Sch:Cyc}
that all the convex polygons can be obtained from regular $p$-gons, 
where $p$ ranges over prime numbers.  We describe these ideas in section
4 and include  proofs 
of Mann's bound on the order of roots of unity satisying a linear equation 
of a given length (see Prop. 4.3) 
and of Schoenberg's result (see Cor. 4.4). 

The results in sections 3 and 4 are used in section 5 
to find the rational planes in an algebraic subset 
$V \subset (\C^*)^r$ (see Prop. 5.2)
and to give formulas for $\deg(p_V)$ 
and bounds on $\per(p_V)$ (see Prop. 5.3 and Theorem 3) in terms of 
defining equations for $V$.  

In section 6 we conclude with some illustrative examples.
For instance, we find $p_V$ when $V$ is the Fermat curve
$$
x^m + y^m = 1, \qquad m \ge 1;
$$
and show (see Example 3)
that, in this case, the bounds for the degree and period
of $p_V$ given
in Theorem 2 are attained.

\heading{Acknowledgement.} I would like to thank E. Jahangard for useful
discussions during the research for this paper.

\section{Notation and main result.}

In this section we will set up some notation to state our main result:
Theorem 2.  This theorem gives bounds on the degree and period of 
$p_V$ rather than exact formulas for them and applies to varieties $V$ 
defined over $\Q$.  A more general result, Theorem 3, is 
stated and proved in section 5, but the bounds in Theorem 2 are 
easier to state and compute.

\heading{Notation.}
Let $\Lambda_r = \C[t_1,\dots,t_r]$ denote the
ring of Laurent polynomials.  A monomial $t_1^{\lambda_1}\dots t_r^{\lambda_r}
\in \Lambda_r$  will be written
as $t^\lambda$ where $\lambda = (\lambda_1,\dots,\lambda_r)$.  Thus,
coordinate functions $t_1,\dots,t_r$ and the usual basis for $\Z^r$
determine a canonical isomorphism between $\Lambda_r$ and $\C[\Z^r]$.

A Laurent polynomial $f \in \Lambda_r = \C[t_1^{\pm 1},\dots,t_r^{\pm 1}]$
will be written as
$$
f(t) = \sum_{\lambda \in \sL(f)} a_\lambda t^{\lambda},
$$
where $\sL(f) \subset \Z^r$ is a finite subset and $a_\lambda \in \C^*$.

For a finite subset $\sF \subset \Lambda_r$, 
our main theorem will describe
the rational planes in the zero set $V(\sF)$ and bounds for 
the degree  and period 
of $p_{V(\sF)}$ in terms of the functions in $\sF$. 
Two main ingredients will be numbers $N[R(\sF)]$ and $D(\sF)$
which can be computed from the number of coefficients of functions 
in $\sF$ and certain subgroups of $\Z^r$ associated to the 
exponents of functions in $\sF$.

Define
$$
R(\sF) = \max\  \{|\sL(f)| : f \in \sF\},
$$
where $|\sL(f)|$ denotes the number of elements in $\sL(f)$.
For any positive integer $R$, let $N[R]$ denote the product of primes
less than or equal to $R$.  

The number $D(\sF)$ takes longer to define.  
For any $f \in \Lambda_r$, 
let $\Pi_f$ be the set of all partitions  $\sP$ of $\sL(f)$
such that each $\nu \in \sP$ has at least two elements
and for any finite subset $\sF = \{f_1,\dots,f_\ell\} \in \Lambda_r$, let
$$
\Pi_{\sF} = \Pi_{f_1} \times \dots \times \Pi_{f_\ell}.
$$
For any partition $\sP \in \Pi_f$, 
let $\Epsilon(\sP,f) \subset \Z^r$ be the subgroup generated by 
$$
\{\lambda-\mu : \ \exists\nu \in \sP,\ \lambda,\mu \in \nu\};
$$
for any $\pi = (\sP_1,\dots,\sP_k) \in \Pi_{\sF}$, let
$\Epsilon(\pi,\sF) \subset \Z^r$ be the sum
$$
\Epsilon(\sP_1,f_1) + \dots + \Epsilon(\sP_k,f_k);
$$
and for any subset $\sU \subset \Pi_{\sF}$, let
$\Epsilon(\sU,\sF) \subset \Z^r$ be the subgroup generated by
$$
\bigcup_{\pi \in \sU}\Epsilon(\pi,\sF).
$$

For any finite abelian group $G$, let $D(G)$
be the largest order of any element in $G$.
Define
$$
D(\sU, \sF) = D(\overline{\Epsilon(\sU,\sF)}/{\Epsilon(\sU,\sF)})
$$
and
$$
D(\sF) = \lcm\  \{D(\sU,\sF) : \sU \subset \Pi_{\sF}\}.
$$

\smallskip

Our main result is the following.

\begin{theorem} Let $V = V(\sF) \subset (\C^*)^r$ be any algebraic subset
defined over $\Q$. 
\begin{description}
\item{(i)}  For each maximal 
rational plane $Q$ in $V$, there is a $\pi \in \Pi_{\sF}$
so that $Q$ is a translate of the affine subtorus of $(\C^*)^r$ defined by 
the set of binomials
$$
\{t^\lambda - 1 :  \lambda \in \overline{\Epsilon(\pi,\sF)} \};
$$
\item{(ii)} 
$$\deg(p_V) \leq r - \min_{\pi\in \Pi_{\sF}} \rank(\overline{\Epsilon(\pi,\sF)}); 
$$
\item{(iii)} 
$$
\per(p_V) \quad | \quad N[R(\sF)]\ D(\sF).
$$
\end{description}
\end{theorem}

This theorem together with finer versions will be proved in section 5.

\section{Character groups and rational planes.}

The main tool we use in this paper is a natural correspondence between
between affine subtori of $(\C^*)^r$ and subgroups of $\Z^r$.  
The relation is explained in section 3.1 and applied in section 3.2
to give properties of rational planes.

\subsection{Subgroups of $\Z^r$ and algebraic subgroups of $(\C^*)^r$.}
For any subgroup $\Epsilon \subset \Z^r$, let $V(\Epsilon)$ be the
closed points of $\Spec (\C[\Z^r/\Epsilon])$.  Then, there is a natural
identification of $V(\Epsilon)$ with the group of characters of
$\Z^r/\Epsilon$, since any closed point corresponds to a ring 
epimorphism $\C[\Z^r/\Epsilon] \rightarrow \C$.  
Define $I_\Epsilon \subset \Lambda_r$ to be the ideal generated by
$$
\{ t^\lambda - 1 : \lambda \in \Epsilon\}.
$$
Then the kernel of the map $\C[\Z^r] \rightarrow \C[\Z^r/\Epsilon]$
is the image of $I_\Epsilon$ in $\C[\Z^r]$ under the identification
$\Lambda_r \iso \C[\Z^r]$.  Thus, there is a natural embedding of $V(\Epsilon)$
into $(\C^*)^r$ with $V(\Epsilon) = V(I_\Epsilon)$ and 
$I_\Epsilon = I(V(\Epsilon))$.  
This embedding preserves group structure as well as the algebraic structure
of $V(\Epsilon)$.

For any algebraic subset $V \subset (\C^*)^r$, define 
$$
\Epsilon(V) = \{\lambda \in \Z^r : t^\lambda - 1\ \mbox{vanishes on}\ V\}.
$$

\begin{lemma} If $\Epsilon \subset \Z^r$ is any
subgroup, then $\Epsilon(V(\Epsilon)) = \Epsilon$.
\end{lemma}

\heading{Proof.}
Observe that for 
any subgroup $\Epsilon \subset \Z^r$, if $t^\lambda-1 \in I_\Epsilon$, 
then $\lambda$ must be a sum of elements in $\Epsilon$.
This gives the inclusion $\Epsilon(V(\Epsilon)) \subset \Epsilon$.
The other inclusion is clear.
\qed

\begin{lemma} Any subtorus $P$ of $(\C^*)^r$ is of the form $V(\Epsilon)$
for some $\Epsilon \subset \Z^r$.
\end{lemma}

\heading{Proof.}
Since $P$ is itself isomorphic to $(\C^*)^s$ for some $0 \leq s \leq r$,
the embedding $\psi : P \rightarrow (\C^*)^r$
defines a ring epimorphism on coordinate rings
$$
\psi^*: \C[\Z^r] \rightarrow \C[\Z^s].
$$
We will show that $\psi^*$ restricted to
$\Z^r$ induces an epimorphism $\Z^r \rightarrow \Z^s$.
This is the same as saying that the components of 
$\psi$ are monomials.  Since $\psi$ preserves multiplication 
it also preserves the multiplication of $(\C^*)$.  Thus, 
the components of $\psi$ must be homogeneous.  Furthermore,
since the components of $\psi$ can have no zeros or poles, other than
the origin, they must be monomials.

Let $\Epsilon \subset \Z^r$ be the kernel of the epimorphism 
$\Z^r \rightarrow \Z^s$ induced
by $\psi$.  Then, the kernel of $\psi^*$ is $I_\Epsilon$, so
$P = V(\Epsilon)$.\qed.

\heading{Remark.}  Another way to say the above is that
there is a naturally duality between algebraic subgroups of $(\C^*)^r$
and quotient groups of $\Z^r$ given by the contravariant functors
$\Hom_{\sA}(-,\C^*)$ and $\Hom(-,\C^*)$, where $\Hom_{\sA}(-,-)$ are
morphisms of algebraic subgroups of $(\C^*)^r$ which preserve both
the algebraic and multiplicative structure.  For brevity and because
they are not necessary for the results of this paper, we omit the details.

\subsection{Rational planes and monomial mappings.}
Given a subgroup $\Epsilon \subset \Z^r$, let 
$\overline{\Epsilon}$ be the subgroup of $\Z^r$ defined by
$$
\overline{\Epsilon} = \{\lambda \in \Z^r : \exists n \in \N,\ 
n\lambda \in \Epsilon \}.
$$

For any subgroup $\Epsilon \subset \Z^r$ 
and $m \in \N$, let
$$
\Epsilon_{m} = \{\lambda \in \Z^r : m\lambda \in \Epsilon \}.
$$
\smallskip

\begin{lemma} For any subgroup $\Epsilon \in \Z^r$, 
$V(\overline\Epsilon)$ is an affine subtorus of $(\C^*)^r$.
Furthermore, there is an injective (non-canonical) endomorphism
$$
T :  \chargp{\OEpsilon} \hookrightarrow (\C^*)^r,
$$
from the character group $\chargp{\OEpsilon}$ into $(\C^*)^r$,
such that $V(\Epsilon)$ decomposes as a disjoint union
$$
V(\Epsilon) = \bigcup_{\eta \in T(\chargp{\OEpsilon})} \eta V(\overline\Epsilon).
$$
\end{lemma}

\heading{Proof.}
The epimorphisms
$$
\Z^r \rightarrow \Z^r/\Epsilon \mapright{} \Z^r/{\overline\Epsilon}
$$
induce inclusions
$$
V(\overline{\Epsilon}) \hookrightarrow V(\Epsilon) \hookrightarrow
(\C^*)^r.
$$
Since $\Z^r/{\overline{\Epsilon}}$ is a free abelian group, 
$V(\overline{\Epsilon})$ is an affine subtorus of $(\C^*)^r$.

The inclusion $\OEpsilon \hookrightarrow \Z^r/\Epsilon$
induces a surjective map on character groups
$$
\psi : V(\Epsilon) \rightarrow  \chargp{\OEpsilon}
$$
whose identity fiber $F_1$ equals $V(\overline\Epsilon)$.  

We need to find a splitting for $\psi$.
Write ${\Z^r}/\Epsilon$ as the product
$$
{\Z^r}/\Epsilon = \Z^s \times G,
$$
where $G$ is finite.  Then $G \iso \OEpsilon$ and 
$\Z^s \iso \Z^r/\overline{\Epsilon}$.
Thus, there is a surjection $\Z^r/\Epsilon \rightarrow \OEpsilon$
whose restriction to $\OEpsilon$ is the identity.
Let $T : \chargp{\OEpsilon}\rightarrow V(\Epsilon)$ be the 
induced endomorphism of character groups.
Then $T$ defines a splitting for $\psi$ and we are done.
\qed
\smallskip

\heading{Notation.}
A map
$$
\psi : \chargp{H} \rightarrow \chargp{G}
$$
between character groups
which is induced by a homomorphism $\psi^* : G \rightarrow H$ 
is called a {\it monomial mapping}, since the induced
map on coordinate rings is just the linear extension of $\psi^*$
to a mapping $\C[G] \rightarrow \C[H]$.

Given a monomial mapping $\psi : (\C^*)^r \rightarrow \chargp{G}$, 
let $\Epsilon(\psi) \subset \Z^r$ be the image of the induced map 
$$
\psi^* : G \rightarrow \Z^r.
$$
\smallskip

\begin{lemma}  Let $\psi : (\C^*)^r \rightarrow \chargp{G}$
be a monomial map and set 
$\Epsilon = \Epsilon(\psi)$.   For each $\mu \in \Tor(\im(\psi))$,
there is a $\tau \in \OEpsilon$
such that, for some $\eta \in (\C^*)^r$ with $\order(\eta) = \order(\mu)\order(\tau)$,
$$
\psi^{-1}(\mu) = \eta V(\Epsilon).
$$
\end{lemma}

\heading{Proof.}   If $m = \order(\mu)$, then the fiber
$\psi^{-1}(\mu)^m$ lies in the identity fiber $V(\Epsilon)$.
Furthermore, if $\eta \in \psi^{-1}(\mu)$ has finite order, then
$m$ divides $\order (\eta)$.  

We have $V(\Epsilon)^m = V(\Epsilon_{m})$, so 
by Lemma 3.3,
$$
V(\Epsilon_{m}) = \bigcup_{\gamma \in T(\overline{\Epsilon}/
{\Epsilon_{m}})}
\gamma V(\overline{\Epsilon}).
$$
Thus, for any $\eta \in \psi^{-1}(\mu)$,
$$
\psi^{-1}(\mu)^m = \eta^m V(\Epsilon)^m = \eta^m V(\Epsilon_{m})
$$
and $\psi^{-1}(\mu)^m$ is a translate of $V(\Epsilon_{m})$
in $V(\Epsilon)$.
Thus, $\psi^{-1}(\mu)^m$ must contain an element of $T(V(\OEpsilon))$.  

This means
we could have chosen $\eta \in \psi^{-1}(\mu)$ such that 
$\eta^m \in T(V(\OEpsilon))$.
If $\Epsilon_{m} = \Epsilon$, then $V(\Epsilon)^m = V(\Epsilon)$
and we could have chosen $\eta$ so that $\eta^m = 1$,
which would imply $\order(\eta) = m$.  Otherwise, set
$\tau = \eta^m$.  We will show 
that  $\order(\eta) = \order(\tau) m$.
We've seen that $m$ divides $\order(\eta)$.  If $k = \order(\eta)/m$,
then $\tau^k = \eta^{mk} = 1$, so $\order(\tau)$ divides $k$.
Since $\eta^{{\ord(\tau)}m} = \tau^{\ord(\tau)} = 1$, we have
$\order(\eta) = \order(\tau)m$.\qed

\begin{corollary} Let $\psi : (\C^*)^r \rightarrow \chargp{G}$ be a 
monomial map and let $\Epsilon = \Epsilon(\psi)$. Then, for
any $\mu \in \im (\psi)$ and connected component
$Q \subset \psi^{-1}(\mu)$, we have 
\begin{description}
\item{(i)} $Q$ is a translate of $V(\overline{\Epsilon})$,
\item{(ii)} $\dim(Q) = r - \rank(\overline{\Epsilon})$, and
\item{(iii)} if $\mu$ has finite order, then 
$$
\ord(\mu)\  | \ \ord (Q) \ |
\ \ord(\mu)D({\overline{\Epsilon}}/\Epsilon),
$$
where $D(\OEpsilon)$ is the largest order of any element of
$\OEpsilon$.
\end{description}
\end{corollary}

We will
now describe the dimensions and orders of  intersections of rational
planes.

\begin{proposition} Let $\Epsilon_1,\dots,\Epsilon_k \subset \Z^r$ be 
any subgroups of $\Z^r$ and let $\eta_1,\dots,\eta_k \in \Tor((\C^*)^r)$.
Let $Q_1,\dots,Q_k \subset (\C^*)^r$ be defined by 
$$
Q_i = \eta_i V(\Epsilon_i), \qquad i=1,\dots,k.
$$
Let $\Epsilon = \Epsilon_1 + \dots + \Epsilon_k$ and
$\eta = (\eta_1,\dots,\eta_k)$.  Let
$$
\begin{array}{rcl}
\rho: \Epsilon_1 \times \dots \times \Epsilon_k &\rightarrow& \Epsilon\\
(\lambda_1,\dots,\lambda_k) &\mapsto& \lambda_1+\dots+\lambda_k
\end{array}
$$
and let $\gamma : \Epsilon_1 \times \dots \times \Epsilon_k 
\rightarrow (\Z^r )^k$ be the inclusion map.
Then $Q_1 \cap \dots \cap Q_k$ is nonempty if and only if 
$\gamma^*(\eta) \in \im(\rho^*)$ and
for any connected component $Q \subset Q_1 \cap \dots \cap Q_k$, 
we have 
\begin{description}
\item{(i)} $Q$ is a translate of $V(\overline{\Epsilon})$,
\item{(ii)} $\dim (Q) = r - \rank (\overline\Epsilon)$, and
\item{(iii)} if $\eta$ has finite order, then
$$
\ord(\eta)\ |\ \ord(Q) \ | \ \order(\eta)D(\OEpsilon).
$$
\end{description}
\end{proposition}

\heading{Proof.}  
Consider the commutative diagram
$$
\cd
{
&\Z^r/{\Epsilon} &\leftarrow &{\Z^r}/{\Epsilon_1} \times\dots\times {\Z^r}/{\Epsilon_k}\cr
&\uparrow{} &&\uparrow{}\cr
&\Z^r &\leftarrow &\Z^r \times \dots \times \Z^r \cr
&\uparrow{}&&\uparrow{}\cr
&\Epsilon &\leftarrow 
&\Epsilon_1 \times \dots \times \Epsilon_r 
}
$$
where the horizontal arrows are given by adding coordinates, 
the bottom vertical arrows are inclusions and the top vertical
arrow are quotient maps.  For $i=1,\dots,k$, let $P_i = V(\Epsilon_i)$.
Then we
have a commutative diagram of induced maps on the associated algebraic sets:
$$
\cd
{
&V(\Epsilon) &\mapright{} &P_1 \times \dots \times P_k \cr
&\downarrow{} &&\downarrow{}\cr
&(\C^*)^r &\mapright{\alpha} &(\C^*)^r\times\dots\times (\C^*)^r\cr 
&\downarrow{\beta} &&\downarrow{\gamma^*}\cr
&\chargp{\Epsilon} &\mapright{\rho^*}&\chargp{\Epsilon_1} 
\times \dots \times \chargp{\Epsilon_k}.}
$$
The preimage $\alpha^{-1}(Q_1\times\dots \times Q_k)$
equals the intersection $Q_1 \cap \dots \cap Q_k$.
If $Q \subset Q_1 \cap \dots \cap Q_k$ is a connected component, then
it is a connected component of a fiber of $\beta$.  The intersection
is nonempty if and only if 
$$
\gamma^*(Q_1\times\dots\times Q_k) \cap \im(\rho^*) \neq \emptyset.
$$
The claim now follows
from Cor. 3.5.\qed
\smallskip

\section{Formal $\Q$-linear relations.}
In this section, we describe $\Q$-linear 
relations between roots of unity and polar rational polygons
following the work of Schoenberg \cite{Sch:Cyc}
and Mann \cite{Mann:LinRels}.   The aim is to give a constructive
method for producing all $\Q$-linear relations and give a 
bound on the orders of roots of unity satisfying a linear equation
in terms of the number of coefficients of the equation.
Other work in this area can be found in {\cite{Len:Van}} and \cite{C-J:Trig}.

For any $n \in N$, let $\zeta_n = \exp(2 \pi \sqrt{-1}/n)$.  For
$z\in \C^*$, let $\theta(z) = \arg(z)/{2 \pi}$.  Thus,
$$
z = |z|\exp(2 \pi \sqrt{-1}\ \theta(z)).
$$

A {\it polar rational polygon}, or prp, is an oriented polygon 
in the complex plane with no parallel edges, such that each edge is a vector
with rational length and  with angle equal to a rational multiple of $2 \pi$.

Then prp's can be put into one-to-one correspondence
with formal
$\Q$-linear relations of the form
\begin{eqnarray}
\sum_{i=1}^k a_i \epsilon_i = 0.
\end{eqnarray}
where $a_i \in \Q^*$ and 
$$
0 \leq \theta(\sign(a_1)\epsilon_1) < \dots
< \theta(\sign(a_k)\epsilon_k) < 1.
$$
Given a prp $T$ defined by (2), we will call $a_i\epsilon_i$ the 
sides of $T$, $|a_i|$ the {\it side lengths} and $\epsilon_i$
the {\it side angles}. 

The {\it order} of $T$, written 
$\ord(T)$, is the 
least $n$ such that, for some root of unity $\eta$, 
we have
$$
\eta^n \epsilon^n = 1,
$$
for all side angles $\epsilon$ of $T$.
The length of $T$, written $\length (T)$ is $k$, the number
of sides.

Any formal $\Q$-linear equation 
$$
\sum_{j=1}^k b_j \eta_j
$$
can be put 
in the form (2) by  the following:  
\begin{description}
\item{(i)} if $\theta(\eta_j)\ (\mod\ 1) \ge 1/2$,
then replace $a_j$ by  $-a_j$ and $\eta_j$ by $\zeta_2 \eta_j$;
\item{(ii)} if $\theta(\eta_j) = \theta(\eta_\ell)$, then replace
$b_j \eta_j + b_\ell \eta_\ell$ by $(b_j + b_\ell)\eta_j$; 
\item{(iii)} remove any summand whose coefficient is zero; and
\item{(iv)} reorder the summands.
\end{description}
It is easy to see that to any formal $\Q$-linear equation there is
a unique corresponding equation of the form (2), which can be obtained
using the above steps, and hence a unique
prp.

It follows that the set of prp's forms a commutative
$\Q[\Tor(\C^*)]$-algebra, where $\Tor(\C^*)$ is the set of 
roots of unity.   One just notes that formal $\Q$-linear equations
are closed under addition, multiplication and scalar multiplication
by elements of $\Tor(\C^*)$.
The algebra operations on a prp will be the same as those for
the corresponding formal $\Q$-linear equation.  The resulting
$\Q$-linear equation may not be of the form (2), but there is
a unique prp associated to it which we take as the result of 
the operation.

The algebra operations have the following geometric interpretations.
Multiplication by a positive rational
$a \in \Q_+$
scales the corresponding polygon by $a$. Multiplication by any root of
unity $\eta \in \Tor(\C^*)$ rotates the polygon by 
the argument of $\eta$.  For example, multiplying the polygon in 
figure 1 by
$\eta = -1$ yields the polygon rotated by 180 degrees shown in figure 2.
$$
\epsffile{neg}
$$
Summing is like taking the union except that one needs to reorder the
sides and get rid of redundancies.

Two prp's will be said to be {\it disjoint} if they do not share any sides
angles.
A prp $T$ is {\it primitive} if there are no disjoint prp's $S$ and $U$ such
that $T = S + U$.
Geometrically, a prp is primitive if there is no
way to rearrange the edges of the polygon to get a union of polygons
joined at vertices.
While any prp has a decomposition into a sum of disjoint primitives, this
decomposition is not necessarily unique.

For example, consider the $\Q$-linear equation  given by expanding out
\begin{eqnarray}
(1  + \zeta_3 + \zeta_3^2)(1 + \zeta_5 + \zeta_5^2 + \zeta_5^3 + \zeta_5^4) = 0.
\end{eqnarray}
The associated prp is given in figure 3.
$$
\epsffile{tot}
$$
The prp defined by (3) can be decomposed into primitive 
prp's as in figure 4 and figure 5.
$$
\epsffile{decomp}
$$
The decomposition in figure 4 comes from writing (3) as
$$
\sum_{i=1}^2 \zeta_3^i\  [1 + \zeta_5 + \zeta_5^2 + \zeta_5^3 + \zeta_5^4] = 0
$$
and figure 5 comes from writing (3) as
$$
\sum_{i=1}^4 \zeta_5^i\ [1 + \zeta_3 + \zeta_3^2] = 0.
$$
\smallskip

The sum of nondisjoint prp's could have smaller length than the
sum of the total as we see in the next example.
Let $A$ be the prp defined by 
$$
\zeta_6 + \zeta_6^2 + (-1) = 0
$$
Then $\zeta_5\ A + T_5$ is the prp (see figure 6)
given by the following $\Q$-linear
equation
$$
1 + (\zeta_6 + \zeta_6^5)\zeta_5 + \zeta_5^2 +  \zeta_5^3 + \zeta_5^4 = 0
$$
or
$$
1 + \zeta_{30} + \zeta_{30}^{11} + \zeta_5^2 + \zeta_5^3 + \zeta_5^4 = 0.
$$
$$
\epsffile{pent}
$$

In general we have the following.

\begin{lemma} The lengths of sums satisfies the following
inequalities
$$
\max\ \{\length(T_1),\length(T_2)\} \leq 
\length(T_1 + T_2) \leq \length (T_1) + \length (T_2),
$$
with equality on the right hand side  if $T_1$ and $T_2$
are disjoint and equality on the left hand side if the set of 
side angles of one prp is contained in the other's.
\end{lemma}

We will make use of the following types of prp's.
For any prime $p$, let $\sigma_p(x)$ be the cyclotomic
polynomial
$$
\sigma_p(x) = 1 + x + \dots + x^{p-1}.
$$
Let $T_p$ be the prp defined by $\sigma_p(\zeta_p) =  0$.

Let $n = mp$, where $p$ is a prime not dividing $m$.
(For example, if $n$ is prime then $m = 1$.)
Then the minimal polynomial $\sigma_{n,p}$ 
for $\zeta_n$ over $\Q[\zeta_m]$ is given by 
$$
\sigma_{n,p}(x) = \zeta_m^{a(p-1)} + \zeta_m^{a(p-2)}x
+ \dots + x^{p-1},
$$
where $a$ is the integer $0 \leq a < m$, such that 
$ap \equiv 1\ (\mod\ m)$.
Let $T_{n,p}$ be the prp defined by $\sigma_{n,p}(\zeta_n) = 0$.
That is, the prp corresponding to the formal $\Q$-linear relation
$$
\sum_{i=0}^{p-1} \zeta_n^{ap(p-1-i) + i} = 0.
$$

\begin{lemma} Let $n$ be an integer, $p$ a prime dividing $n$, such
that $p^2$ doesn't divide $n$.  The prp $T_{n,p}$ is a multiple of $T_p$
by some element of $\Q[\Tor(\C^*)]$.
\end{lemma}

\heading{Proof.}  Let $m = n/p$ and let $a$ be such that $ap \equiv 1
\ (\mod\ m)$.  Let $r$ be any integer such that $ap = 1 + mr$.
Then, since $\zeta_n^m = \zeta_p$, we have 
\begin{eqnarray*}
\zeta_n^{ap(p-1-i) + i} &=& \zeta_n^{mr(p-1-i) + (p-1)}\\
&=& \zeta_p^{r(p-1-i)} \zeta_n^{p-1}.
\end{eqnarray*}
Thus, $T_{n,p}$ is defined by the formal linear relation
$$
\zeta_n^{p-1}\sum_{i=0}^{p-1} \zeta_p^{r(p-1-i)} = 0.
$$
Since $p$ doesn't divide $r$, $r(p-1-i)$ ranges in $0,\dots,p-1$
as $i$ ranges in $0,\dots,p-1$.  Therefore, 
$$
T_{n,p} = \zeta_n^{p-1} T_p.
$$
\qed

Our proof of the following result uses essentially the same ideas
as Mann uses in (\cite{Mann:LinRels}, Theorem 1), except that we 
take more advantage of the structure of prp's as a 
$\Q[\Tor(\C^*)]$-algebra.

\begin{proposition} Let $T$ be a primitive prp of length $r$
and order $n$. Then we have
\begin{description}
\item{(i)} $n$ is square free;
\item{(ii)} for any prime $p$ dividing $n$, $p \leq r$;
\item{(iii)} for any prime $p$ dividing $n$, 
$T$ is a $\Q[\Tor(\C^*)]$-linear combination
of $T_{n,p}$ and prp's of order $n/p$.
\end{description}
\end{proposition}

\heading{Proof.}
Let
\begin{eqnarray}
\sum_{i=1}^r a_i \epsilon_i = 0
\end{eqnarray}
be the $\Q$-linear relation corresponding to $T$ of the form (2).
If $T$ has order $n$, then by 
multiplying $T$ by an element of $\Tor(\C^*)$, if necessary, we
can assume that all the $\epsilon_i$ in (4) have order $n$.
We will show that $n$ is square free and for any prime $p$
dividing $n$, $p \leq r$.

Let $p$ be a prime dividing $n$ and let $m = n/p$.  Since 
$\zeta_n^p = \zeta_m$, we can write 
$$
\epsilon_i = \eta_i\zeta_{n}^{\alpha_i},
$$
for $i=1,\dots,r$, where $\eta_i \in \Q[\zeta_{m}]$ and
$0 \leq \alpha_i \leq p-1$.  Define
$$
q_{T,n,p}(x) = \sum_{i=1}^r a_i \eta_i x^{\alpha_i}.
$$
This is a polynomial with coefficients in $\Q[\zeta_{m}]$
satisfied by $\zeta_{n}$ and hence the minimal polynomial
$\sigma_{n,p}(x)$ for $\zeta_{n}$ over $\Q[\zeta_{m}]$
divides $q_{T,n,p}(x)$.
If $p^2$ divides $n$ then $\deg(\sigma_{n,p}) = p$ which is
strictly greater than the degree of $q_{T,n,p}(x)$, and hence 
$q_{T,n,p}(x)$ is identically zero.  Since $T$ is
primitive, this means that all the $\alpha_i$ are the
same, but then we can multiply $T$ by $\zeta_p^{-1}$ to get
a prp of order $n/p$ which is a contradiction, since the order
of a prp is preserved under multiplication by roots of unity.

Thus, $n$ is square free (proving (i)) and $\sigma_{n,p}(x)$ divides
$q_{T,n,p}(x)$.  
Set 
$$
A_\alpha = \sum_{\alpha_i = \alpha} a_i \eta_i 
$$
and write $q_{T,n,p}(x)$ as 
$$
q_{T,n,p}(x) = \sum_{\alpha = 0}^{p-1} A_\alpha x^\alpha.
$$
Since $q_{T,n,p}(x)$ is not identically
zero and has degree less than or equal to $p-1$, 
we have
\begin{eqnarray}
q_{T,n,p}(x) = B \sigma_{n,p}(x)
\end{eqnarray}
for some invertible $B \in \Q[\zeta_m]$.

By (5), evaluating $B^{-1}\ A_\alpha$ as a 
complex number gives the $\ell$th
coefficient of $\sigma_{n,p}(x)$.  
Thus,  the formal $\Q$-linear equation
$$
B^{-1}\ A_\alpha -  \zeta_m^{a(p-1-\ell)} = 0 
$$
defines a prp, which we'll call $T_\alpha$. 
Then $T$ has the following decomposition into a sum of prp's 
$$
T = B\ \sum_{\alpha = 0}^{p-1} T_\alpha + B\ T_{n,p}.
$$
Since $\order(T_\alpha) = m$, we have proved (iii).

By Lemma 4.1, $\length (T) \ge \length (T_{n,p}) = p$, so $p \leq r$
which proves (ii).
\qed

Proposition 4.3 leads to a geometric method for generating all prps as we see
in the following result, originally proved by Schoenberg in \cite{Sch:Cyc}.

\begin{corollary} The set of prp's, considered as a $\Q[\Tor(\C^*)]$-module,
is generated by the set of $T_p$, where $p$ ranges over all primes.
\end{corollary}

\heading{Proof.} By Lemma 4.2, 
$T_{n,p}$ is a multiple of $T_p$ by an element of $\Q[\zeta_{n/p}]$. 
The rest follows from Prop. 4.3 and induction on the number
of primes dividing $n$.
\qed

It follows that geometrically one can construct all prp's by starting
with regular $p$-gons, for primes $p$, and doing the operations of
rotation, stretching and ``adding", where ``adding"  means taking the
union of sides, grouping like side angles, getting rid of sides of
zero length and reordering the sides, as described above.

\section{Behavior of torsion points.}
We start in this section with a review of the proof of Theorem 1.  
Then, using the
relation between finitely generated groups and rational planes
developed in section 3, we study the degree and period of $p_V$.
Finally, we give a proof of Theorem 2.

\subsection{Proof of Theorem 1.}

Let $V = V(f_1,\dots,f_r) \subset (\C^*)^r$ be any algebraic subset.
We want to find a finite number of rational planes $Q_1,\dots,Q_\ell$
contained in $V$ so that 
$$
\Zar{Tor(V)} = \bigcup_{i=1}^\ell  Q_i.
$$

\heading{Step 1.}  We can reduce to the case of hypersurfaces, since 
$V = V(f_1) \cap \dots \cap V(f_r)$ and, by Prop. 3.6, a
finite intersection of rational planes is a finite union of rational
planes.   
\smallskip

\heading{Step 2.}  Let $V = V(f)$, where $f \in \Lambda_r$ defined by
$$
f = \sum_{\lambda \in \sL(f)} a_\lambda t^\lambda, \quad a_\lambda \in \C^*,
$$
for $\sL(f) \subset \Z^r$ is a finite subset.  We'll show that
$$
\Zar{\Tor(V)} = \bigcup_{i=1}^k \psi_i^{-1}(S_i),
$$
for some finite collection of monomial maps $\psi_i: (\C^*)^r \rightarrow
\chargp{G_i}$, some character groups $\chargp{G_i}$, and some finite sets
of torsion elements $S_i \subset \Tor(\chargp{G_i})$. 
It then follows from Cor. 3.5 that $\Zar{\Tor(V)}$ is a finite union of 
rational planes.

Let $\sL = \sL(f)$ and
let $\ell(f)$ be the linear polynomial
$$
\ell(f) = \sum_{\lambda \in \sL} a_\lambda x_\lambda
$$
where $x_\lambda$
are independent variables, for $\lambda$ ranging in
$\sL$.
Then 
$$
V(f) = \alpha_f^{-1}(V(\ell(f))), 
$$
where $\alpha_f : (\C^*)^r \rightarrow (\C^*)^{\sL}$
is the monomial map induced by
\begin{eqnarray*}
\Z^{\sL} &\rightarrow& \Z^r\\
e_\lambda &\mapsto& \lambda.
\end{eqnarray*}
\smallskip

\heading{Step 3.}  Let $\Pi_f$ be the set of 
partitions of $\sL$, where for all $\nu \in \sP \in \Pi_f$,
$\nu$ has at least 2 elements.
Fix $\sP \in \Pi_f$.
Let $L_{\sP}$ be the system of linear equations
$$
\ell_\nu(f) = \sum_{\lambda \in \nu}a_{\lambda}t_{\lambda},
$$
for each $\nu \in \sP$. 
This defines a system of linear equations
defined on
$$
W(\sP)  = \prod_{\nu \in \sP} (\C^*)^{\nu}.
$$
Let $\sT(\sP,f)$ be the set of 
elements $(\epsilon_\lambda)_{\lambda \in \nu, \nu \in \sP} \in W(\sP)$ 
satisfying the system $L_{\sP}$: 
$$
\sum_{\lambda \in \nu} a_\lambda \Epsilon_\lambda = 0,
\qquad\mathrm{for all}\ \nu \in \sP.  
$$
Then 
$$
V(\ell(f)) = \bigcup_{\sP \in \Pi_f} \sT(\sP,f).
$$

For each $\nu \in \sP$, choose $\lambda_\nu \in \nu$
and let $\nu^* = \nu \setminus \{\lambda_\nu\}$.
Define
$$
W(\sP)^* = \prod_{\nu \in \sP} (\C^*)^{\nu^*}.
$$
Let
$$
\beta_\nu : (\C^*)^\nu \rightarrow (\C^*)^{\nu^*}
$$
be the map defined by
\begin{eqnarray*}
\Z^{\nu^*} &\rightarrow& \Z^\nu\\
e_\lambda &\mapsto& e_\lambda - e_{\lambda_\nu}.
\end{eqnarray*}
Let 
$$
\beta_\sP : W(\sP) \rightarrow W(\sP)^*
$$
be the product of the maps $\beta_{\nu}$.
Let $\sV(\sP,f)$ be the set of solutions 
$(\Epsilon_\lambda)_{\lambda \in \nu^*, \nu \in \sP} \in W(\sP)^*$ 
to the inhomogeneous equations
$$
\sum_{\lambda \in \nu^*} (-a_{\lambda}/{a_{\lambda_\nu}}) 
t_\lambda  = 1, \qquad \nu \in \sP.
$$
Then $\sT(\sP,f) = \beta_{\sP}^{-1}(\sV(\sP,f))$, so we have
\begin{eqnarray}
V(f) = \bigcup_{\sP\in \Pi_f} \alpha_f^{-1} \beta_{\sP}^{-1}(\sV(\sP,f)).
\end{eqnarray}
\smallskip

\heading{Step 4.} Let
$$
\psi_{\sP,f} : (\C^*)^r \rightarrow W(\sP)^*
$$
be the composition $\psi_{\sP,f} = \beta_{\sP}\circ \alpha_f$.
Then $\psi_{\sP,f}$ is induced by the maps
\begin{eqnarray*}
\bigoplus_{\nu \in \sP} \Z^{\nu^*} &\rightarrow& \Z^r\\
e_\lambda &\mapsto& \lambda - \lambda_\nu \quad \mbox{if $\lambda \in \nu$.}
\end{eqnarray*}
Equation (6) now becomes
\begin{eqnarray}
V(f) = \bigcup_{\sP\in \Pi_f}\psi_{\sP,f}^{-1}(\sV(\sP,f)).
\end{eqnarray}
Since monomial maps take torsion points to torsion points, (7) implies
$$
\Zar{\Tor(V(f))} = \bigcup_{\sP\in \Pi_f} \psi_{\sP,f}^{-1}(\Tor(\sV(\sP,f)))
$$

If $\sP$ and $\sP'$ are two elements in $\Pi_f$, 
$\sP'$ is a {\it refinement} of $\sP$  if 
for any $\nu' \in \sP'$,
there is a $\nu \in \sP$ such that $\nu' \subset \nu$.
A refinement $\sP'$ of $\sP$ is {\it proper} if $\sP'$ is not equal to $\sP$.
If $\sP'$ is a refinement of $\sP$, then 
$$
\sV(\sP',f) \subset \sV(\sP,f).
$$
Let $\sV_m(\sP,f)$ be the set of elements of $\sV(\sP,f)$ which are
not in $\sV(\sP',f)$ for any proper refinement $\sP'$ of $\sP$.
Let $S(\sP,f) = \Tor(\sV_m(\sP,f))$.
Then we have
\begin{eqnarray}
\Zar{\Tor(V(f))}
= \bigcup_{\sP \in \Pi_f} \psi_{\sP,f}^{-1}(S(\sP,f)).
\end{eqnarray}
Thus, the proof of Theorem 1 is completed by the following lemma.

\begin{lemma} (Laurent, Sarnak)  The set $S(\sP,f)$ is  
finite.
\end{lemma}

This follows from (\cite{Laur:Equ}, Theorem 1) and (\cite{A-S:Betti},
Lemma 3.1). \qed

The rational planes inside a given algebraic subset of $(\C^*)^r$ can now be 
described in the following way.

\begin{proposition}
Let $V\subset (\C^*)^r$ be any algebraic subset defined
by Laurent polynomials $\sF = \{f_1,\dots,f_k\} \subset \Lambda_r$.
Let $\Pi_\sF =\Pi_{f_1} \times \dots \times \Pi_{f_k}$.
Then
\begin{eqnarray}
\Zar{\Tor(V)}
= \bigcup_{\sP_1 \times \dots \times \sP_k \in \Pi_{\sF}}\  \bigcap_{i=1}^k
\psi_{\sP_i,f_i}^{-1}(S(\sP_i,f)).
\end{eqnarray}
\end{proposition}

\heading{Proof.}
We have
$$
\Zar{\Tor(V)} = \bigcap_{i=1}^k \Zar{\Tor(V(f_i))}
$$
so the statement follows from (8).
\qed
\smallskip

\subsection{Degree and Periodicity.}

We will now give some results  concerning the degree and period of
$p_V$. 

First let $f \in \Lambda_r$ be a Laurent polynomial.  We will begin
by studying the hypersurface $V(f) \subset (\C^*)^r$ defined by $f$.
For a partition $\sP \in \Pi_f$, define
$$
\delta_{\sP,f} = \left \{ \begin{array}{ll} 
1 &\mbox{if $S(\sP,f) \neq \emptyset$,}\\
0 &\mbox{otherwise.}
\end{array}\right .
$$
For $\sP \in \Pi_f$, let 
$\Epsilon(\sP,f)  \subset \Z^r$ be the subset generated
by 
$$
\{ \lambda - \mu : \exists \nu \in \sP, \lambda,\mu \in \nu\}
$$
Then $\Epsilon(\sP,f) = \Epsilon(\psi_{\sP,f})$, where $\psi_{\sP,f}$
is the map defined in the previous section,
so by Cor. 3.5, any connected component of a fiber of $\psi_{\sP,f}$
is a translate of $V(\overline{\Epsilon(\sP,f)})$.

For any subset $\sU \subset \Pi_f$, let
$\Epsilon(\sU,f) \subset \Z^r$ be the subgroup generated by
$$
\bigcup_{\sP \in \sU} \Epsilon(\sP,f).
$$
Then by  Prop. 3.6 any rational plane $Q$ in the intersection
$$
\bigcap_{\sP \in \sU} F_{\sP} 
$$
where each $F_{\sP}$ is a fiber of $\psi_{\sP,f}$, is a translate
of $V(\overline{\Epsilon(\sU,f)})$.

Let
$$
D(\sU,f) = D(\overline{\Epsilon(\sU,f)}/{\Epsilon(\sU,f)})
$$
and let
$$
D(f) = \lcm \ \{D(\sU,f) : \sU \subset \Pi_f\}.
$$
Let
$$
M(\sP,f) = \max \ (\{\order(x) : x \in  S(\sP,f)\} \cup \{1\})
$$
and let
$$
M(f) = \lcm_{\sP \in \Pi_f}\  M(\sP,f).
$$

\begin{proposition}  Let $f \in \Lambda_r$ any Laurent polynomial
and $V = V(f)$.  Then 
any rational plane $Q \subset V(f)$
is a translate of $V(\overline{\Epsilon(\sP,f)})$, for some
partition $\sP \in \Pi_f$.
Furthermore,
$$
\deg(p_{V(f)}) = \max_{\sP \in \Pi_f}\ (r - \rank(\overline{\Epsilon(\sP,f)})
\delta_{\sP,f}
$$
and
$$
\per(p_{V(f)}) \quad | \quad M(f)D(f).
$$
\end{proposition}

\heading{Proof.} Recall that
$$
\Zar{\Tor(V(f))} = \bigcup_{\sP\in \Pi_f} \psi_{\sP,f}^{-1}(S(\sP,f)).
$$
Thus, the degree of $p_{V(f)}$ is the maximum dimension of rational
planes in $\psi_{\sP,f}^{-1}(S(\sP,f))$.  
Since $\Epsilon(\psi_{\sP,f}) = \Epsilon(\sP,f)$,
by Cor. 3.5,  if $S(\sP,f)$ is not empty, then 
for any rational plane $Q \subset \psi_{\sP,f}^{-1}(S(\sP,f))$,
we have
$$
\dim(Q) = r - \rank(\overline{\Epsilon(\sP,f)}).
$$
This gives the formula for $\deg(p_{V(f)})$.
The period of $p_V$ depends on the orders of rational planes
in intersections of fibers of $\psi_{\sP,f}^{-1}$ over
$S(\sP,f)$.  Any such intersection is empty if it involves
more than one fiber of $\psi_{\sP,f}$, for some $\sP$.
Let $\sU \subset \Pi_f$ be any subset such that, for any
$\sP \in \sU$, $\delta_{\sP,f} = 1$.  Choose a connected
component $Q_\sP \subset \psi_{\sP,f}^{-1}(S(\sP,f))$  
for each $\sP \in\sU$.  
Let
$$
M = \lcm_{\sP \in \sU} \ \order(\psi_{\sP}(Q_{\sP})).
$$
Then, by Prop. 3.6, for any rational plane $Q$ in  the intersection 
$$
\bigcap_{\sP \in \sU}Q_{\sP},
$$
we have
$$
\ord(Q) \ | \ M\ D(\sU,f).
$$
Taking the least common multiple of both sides, the claim follows.
\qed
\smallskip

While the various groups denoted by $\Epsilon(-)$ can be computed
routinely, the sets $S(\sP,f)$ need to be studied in a case by
case manner.  When $f$ is defined over $\Q$ the results described in
section 4 are aids to computation.

In particular, we have the following result of Mann \cite{Mann:LinRels}
(cf. Prop. 4.3).

\begin{proposition} (Mann)
Let $\epsilon = (\epsilon_1,\dots,\epsilon_r)$ be any
finite order maximal solution to a linear equation
\begin{eqnarray*}
\sum_{i=1}^r a_i \epsilon_i = 1,
\end{eqnarray*}
with $a_i \in \Q$ for all $i=1,\dots,r$.  Let $n$ be the order
of $\epsilon$ as an element of $(\C^*)^r$.  Then $n$ is square
free and, if $n=p_1\dots p_k$ is a factorization, then
$p_i \leq r+1$ for $i=1,\dots,k$.
\end{proposition}

As a consequence of Prop. 5.4, we have the following bound on $M(\sP,f)$.
For any partition $\sP$, let
$$
R(\sP) = \max\ \{|\nu | : \nu \in \sP\}.
$$
For any positive integer $R$, let $N[R]$ be the product of distinct
primes less than or equal to $R$.

\begin{corollary} If $f \in \Lambda_r$ is defined over $\Q$ and has 
$R$ coefficients, then 
$$
M(\sP,f) \ | \  N[R(\sP)]
$$
for any partition $\sP \in \Pi_f$.
\end{corollary}

We will now give bounds on $p_{V(\sF)}$, where $\sF = \{f_1,\dots,f_k\}
\subset \Lambda_r$ is
any finite subset.

For any $\pi= (\sP_1,\dots,\sP_k) \in \Pi_{\sF}$, let
\begin{eqnarray*}
\Epsilon(\pi,\sF) &=& \Epsilon(\sP_1,f_1) + \dots + \Epsilon(\sP_k,f_k),\\
S(\pi,\sF) &=& S(\sP_1,f_1) \times \dots \times S(\sP_k,f_k),\\
R(\pi)  &=& \max\ \{R(\sP_i) : i=1,\dots,k\}, \\
M(\pi,\sF) &=& \lcm \ \{\ord(x) : x \in S(\pi,\sF)\}, \ \mbox{and}\\
M(\sF) &=& \lcm\  \{M(\pi,\sF) : \pi \in \Pi_{\sF}\}.
\end{eqnarray*}

For any subset $\sU \subset \Pi_\sF$, define
$\Epsilon(\sU,\sF) \subset \Z^r$ to be the subgroup generated by
$$
\bigcup_{\pi \in \sU}\Epsilon(\pi,\sF).
$$
Let
\begin{eqnarray*}
D(\sU, \sF) &=& D(\overline{\Epsilon(\sU,\sF)}/{\Epsilon(\sU,\sF)})\\
D(\sF) &=& \lcm\  \{D(\sU,\sF) : \sU \subset \Pi_{\sF}\}.
\end{eqnarray*}
Thus, for any $\lambda \in \overline{\Epsilon(\sU,\sF)}/{\Epsilon(\sU,\sF)}$,
we have $\order(\lambda)\ | \ D(\sF)$.

Let $\rho$  be defined by
\begin{eqnarray*}
\rho: \Epsilon(\sP_1,f_1) \times \dots \times \Epsilon(\sP_k,f_k)
&\rightarrow &
\Epsilon(\sP_1,f_1) + \dots + \Epsilon(\sP_k,f_k)\\
(\lambda_1,\dots,\lambda_k)&\mapsto& \lambda_1 + \dots + \lambda_k
\end{eqnarray*}
and let 
$$
\gamma: \Epsilon(\sP_1,f_1) \times \dots \times \Epsilon(\sP_k,f_k)
\hookrightarrow
Z^r \times \dots \times \Z^r
$$
be the inclusion map.  Let
$$
\delta_{\pi,\sF} = \left\{\begin{array}{ll} 1 &\quad\mbox{if 
$\gamma^*(S(\pi,\sF)) \cap \im (\rho^*) \neq \emptyset$,} \\
0 &\quad\mbox{otherwise.}
\end{array}\right .
$$

Prop. 5.3 extends to arbitrary algebraic subsets as follows.

\begin{theorem}  Let 
$\sF \subset \Lambda_r$ be a finite subset let and
$V = V(\sF)$.  Then any rational plane
$Q \subset V$ is a translate of $V(\overline{\Epsilon(\pi,\sF)})$ for 
some partition $\pi \in \Pi_{\sF}$.
Furthermore,
$$
\deg(p_V) = \max_{\pi \in \Pi_{\sF}}
\ (r - \rank(\overline{\Epsilon(\pi,\sF)})) \delta_{\sP,\sF}
$$
and
$$
\per(p_V) \ | \ M(\sF)\ D(\sF).
$$
\end{theorem}

\heading{Proof.}  From (9), we know that $\Zar{\Tor(V(\sF))}$ 
is the union over $\pi \in \Pi_{\sF}$ of 
$$
Q(\pi,\sF) = \bigcap_{i=1}^k \psi_{\sP_i,f}^{-1}(S(\sP_i,f)).
$$
Hence $\deg(p_V)$ is the maximum dimension of any rational plane
$Q \subset Q(\pi,\sF)$.
By Prop. 3.6, $Q(\pi,\sF)$ is nonempty if and only if 
$\delta_{\pi,\sF} = 1$. 
Also by Prop. 3.6, any rational plane $Q \subset Q(\pi,\sF)$
is a translate of $V(\overline{\Epsilon(\pi,\sF)})$ and its dimension
equals $r-\rank(\overline{\Epsilon(\pi,\sF)})$. This gives the formula
for the degree of $p_V$.

The period $\per(p_V)$ is the $\lcm$ of $\ord(Q)$, where $Q$ is
a rational plane in the intersection $\bigcap_{\pi \in \sU} Q_{\pi}$,
where $Q_{\pi} \subset Q(\pi,\sF)$,
is some choice of rational planes
and $\pi$ ranges
in a subset $\sU \subset \Pi_\sF$.
Note that 
$$
M(\sU,\sF) = \lcm\ \{ M(\pi,\sF) : \pi\in \sU\}
$$
divides $M(\sF)$.
By Prop. 3.6, we have
$$
\ord(Q) \ | \ M(\sU,\sF)\ D(\sU,\sF) 
$$ 
and taking the least common multiple of both sides, we get 
the bound for the period of $p_{V(f)}$.
\qed

\subsection{Proof of Theorem 2.}  Theorem 2 follows easily from
Theorem 3 and Corollary 5.5.   To prove (i), we need to find
find polynomials defining the rational planes in $V$.  Since by
Theorem 3, any rational plane $Q$ in $V$ is of the form $Q = \eta 
V(\Epsilon)$,
where $\Epsilon = \overline{\Epsilon(\pi,\sF)}$ and $\eta$ is a finite
order element of $(\C^*)^r$.  By Lemma 3.1, $V(\Epsilon)$ is
defined by $I_\Epsilon$ and (i) follows.
The bound for $\deg(p_V)$ in (ii) 
comes from ignoring the $\delta_{\sP,V}$ in Theorem 3.  
By Corollary 5.5, we have the inequality
$$
M(\sF) \ | \   N[R(\sF)],
$$
where $R(\sF)$ is the maximum number of coefficients of a Laurent polynomial
$f \in \sF$.  This implies the bound for $\per(p_V)$ in (iii). 
\qed

\section{Examples.}
In this section we consider some algebraic subsets $V \subset (\C^*)^r$
defined by a finite set of Laurent polynomials $\sF \subset \Lambda_r$ 
and study the torsion points on $V$ in terms of $\sF$.  
We start with some notation.
For any $m \in \N$, let  $\sC_m : \N \rightarrow \N$ 
be the periodic function defined by
$$
\sC_m(n) = \left\{\begin{array}{ll}1 &\quad\mbox{if $m$ divides $n$,}\\
0 &\quad\mbox{otherwise}\end{array}\right .
$$

For any root of unity $\epsilon$, let $\theta(\epsilon) \in \Q/\Z$
be such that
$$
\epsilon =  \exp(2\pi\ \sqrt {-1} \ \theta(\epsilon)).
$$
For any positive integer $n$, denote by 
$\zeta_{n}$ the primitive $n$th root of unity
such that $\theta(\zeta_n) = 1/n$.
For any $\theta \in (\Q/\Z)^r$, let $\order(\theta)$ be the
least positive integer $n$, such that $n\theta = 0 (\mod\ 1)$.

For any torsion point $\epsilon = (\epsilon_1,\dots,\epsilon_r)$
in $(\C^*)^r$ define
$$
\theta({\epsilon}) = (\theta(\epsilon_1),\dots,\theta(\epsilon_r))
\ \in \ (\Q/\Z)^r.
$$
Then
$$
\order(\epsilon) = \order(\theta(\epsilon)).
$$

A binomial equation
\begin{eqnarray}
t^\lambda - \zeta = 0
\end{eqnarray}
where 
$\lambda=(\lambda_1,\dots,\lambda_r) \in \Z^r$ and $\zeta$ is a root
of unity $\zeta = \zeta_{n}^k$, $0 \leq k \leq n-1$,
corresponds to a linear equation
\begin{eqnarray}
\lambda_1 \theta_1 + \dots + \lambda_r \theta_r = k/n\ (\mod\ 1),
\end{eqnarray}
in the sense that $\epsilon$ satisfies (10) if and only if $\theta(\epsilon)$ is
a solution to (11).
We'll call equation (11) the {\it exponential form} of equation (10).

We'll use the following simplification of  Theorem 3
in Examples 2 and 3.

\begin{proposition}  Suppose $\Pi_{\sF}$ contains only one element
$\pi$.  
Let 
$$
M(\sF) = \lcm \ \{\order(x) : x \in S(\pi,\sF)\} \cup \{1\}
$$ 
and let
$$
D(\sF) = D(\overline{\Epsilon(\pi,\sF)}/{\Epsilon(\pi,\sF)}).
$$
If $S(\pi,\sF) = 0$, then $p_{V(\sF)} = 0$.
Otherwise,
$$
\deg (p_{V(\sF)}) = r - \rank(\overline{\Epsilon(\pi,\sF)}).
$$
and 
$$
M(\sF) \ | \ \order(p_{V(\sF)}) \ |\ M(\sF)\ D(\sF).
$$
\end{proposition}

\heading{Proof.}  We have
$$
\Zar{\Tor(V)} = \psi_{\pi,\sF}^{-1}(S(\pi,\sF)),
$$
which is a union of disjoint rational planes
contained in fibers over the points in $S(\pi,\sF)$.
Since $\Epsilon(\pi,\sF) = \Epsilon(\psi_{\pi,\sF})$,
the rest follows from Cor. 3.5.
\qed

\heading{Example 1.}   We begin with examples 
where the degree bound
in Theorem 2 is attained.  Let $r$ be an even number $r = 2k$.
Let
$$
f = \sum_{i=1}^r (-1)^i t_i \in \Lambda_r
$$
and $V = V(f)$.
\smallskip

We can see immediately that $V(f)$ contains the affine subtorus 
$$
P = V(\{t_it_{i+1}^{-1} - 1 : i=2j-1, j=1,\dots,k\}).
$$
The dimension of $P$ equals the dimension of solutions 
$\theta=(\theta_1,\dots,\theta_r) \in (\Q/\Z)^r$ such that
$$
\theta_{2j-1} - \theta_{2j} = 0\ (\mod\ 1), \qquad \mbox{for all
$j = 1,\dots,k$}.
$$
The dimension is clearly $r - k = k$, so $\deg (p_{V(f)}) \ge k$. 

Now take any partition $\sP \in \Pi_f$.  Then 
$\Epsilon(\sP,f) = \overline{\Epsilon(\sP,f)}$ and
$$
\rank(\Epsilon(\sP,f)) = \sum_{\nu \in \sP} (|\nu|-1) \ge k.
$$
Therefore, 
$$
\deg (p_{V(f)}) \leq  n - k = k
$$
and we have equality.
\smallskip

\heading{Example 2.}  Let $V(f) \in (\C^*)^r$ be defined by
$$
f = a_1 t^{\lambda_1} + a_2 t^{\lambda_2}
+ a_3 t^{\lambda_3}, \qquad a_i \in \Q^*
$$
Since there are only three coefficients there is only one parition
in $\Pi_f$, namely $\sP = \{1,2,3\}$.

We have
\begin{eqnarray*}
&&\Epsilon = \Epsilon(\sP,f) 
= \{\lambda_1 - \lambda_3, \lambda_2 - \lambda_3\}\\
&& L_{\sP,f} : \frac{-a_1}{a_3}x + \frac{-a_1}{a_2}y - 1\\
\end{eqnarray*}

The following is a table of all possible
values for $a=-a_1/{a_3}$ and $b=-a_2/a_3$ which give nonempty
$S(\sP,f)$.

\bigskip
\begin{center}
\begin{tabular}{|l|l|}\hline
$(a,b)$ & $S(\sP,f)$\\
\hline
$(1,1)$&$(\zeta_{6},\zeta_{6}^{-1}),(\zeta_{6}^{-1},\zeta_{6})$\\
$(1,-1)$&$(\zeta_{6},\zeta_{3}),(\zeta_{6}^{-1},\zeta_{3}^{-1})$\\
$(-1,1)$&$(\zeta_{3},\zeta_{6}),(\zeta_{3}^{-1},\zeta_{6}^{-1})$\\
$(-1,-1)$&$(\zeta_{3},\zeta_{3}^{-1}),(\zeta_{3}^{-1},\zeta_{3})$\\
$(1/2,1/2)$&$(1,1)$\\
$(1/2,-1/2)$&$(1,-1)$\\
$(-1/2,1/2)$&$(-1,1)$\\
$(-1/2,-1/2)$&$(-1,-1)$\\
\hline
\end{tabular}
\end{center}

\bigskip

Thus, for example, if $f = t_1 + t_2 + t_3$, then 
$$
p_{V(f)}(n) = 2n \ \sC_3(n).
$$

If $\sF \subset \Lambda_r$ is any finite collection
of Laurent polynomials with three coefficients, then $\Pi_{\sF}$ contains
a single element $\pi$ so we can use Prop. 6.1.

Let $V = V(\sF)$, $\sF = \{f_1,f_2\} \subset 
\Lambda_r$, where
\begin{eqnarray*}
f_1 &=& t_1t_3 + t_4 + \alpha\\
f_2 &=& t_1 + t_2 + \beta t_3
\end{eqnarray*}
where $\alpha,\beta \in \{ \pm 1\}$.

Then $\Pi_{\sF}$ has a single element $\pi = (\sP,\sP)$, where $\sP  
= \{1,2,3\}$. Since
\begin{eqnarray*}
&&\Epsilon = \Epsilon(\pi,\sF) = \Epsilon(\sP,f_1) + \Epsilon(\sP,f_2)\\
&&\quad = \{(1,0,1,0),(0,0,0,1), (1,0,-1,0),(0,1,-1,0)\}\\
&&\overline{\Epsilon} = \Z^4,\\
&&D(\sF) = D(\overline{\Epsilon}/\Epsilon) = D(\Z/{2\Z}) = 2,
\end{eqnarray*}
by Prop. 6.1,
the bounds on the degree and period of $p_V$ are
$$
0 \leq \deg(p_V) \leq r - \rank(\overline{\Epsilon})=4 - 4 = 0\\
$$
which implies $\deg(p_V) = 0$ and
$$
M(\sF) \ | \ \per(p_V) \ | \ 2M(\sF).
$$

For any $\alpha$ and $\beta$, the exponential form of equations
defining
$$
\Zar{\Tor(V)} = \psi_{\pi,\sF}^{-1}(S(\pi,\sF)),
$$
are the linear equations
\begin{eqnarray*}
\theta_1 + \theta_3 &=& c\ (\mod\ 1)\\
\theta_4 &=& -c\ (\mod\ 1)\\
\theta_1 - \theta_3 &=& d\ (\mod\ 1)\\
\theta_2 - \theta_3 &=& -d\ (\mod\ 1).\\
\end{eqnarray*}
where $c$ and $d$ range in $A \times B$, and depend on $\alpha$ and
$\beta$.  The following table shows the $p_V$ corresponding to the
different choices of $\alpha$ and $\beta$.

\begin{center}
\begin{tabular}{|l|l|l|l|l|l|}\hline
$(\alpha,\beta)$ & $A$ & $B$ & $M(\sF)$ & $p_V$ & $\per(p_V)$\\
\hline
$(1,1)$ & $\{1/3,2/3\}$ & $\{1/3,2/3\}$ & $3$ & $4\ \sC_3 + 4\ \sC_6$
&$6$\\
$(1,-1)$ & $\{1/3,2/3\}$ & $\{1/6,5/6\}$ & $6$ & $8\ \sC_{12}$  
&$12$\\
$(-1,1)$ & $\{1/6,5/6\}$ & $\{1/3,2/3\}$ & $6$ & $8\ \sC_{12}$
&$12$\\
$(-1,-1)$ & $\{1/6,5/6\}$ & $\{1/6,5/6\}$ & $6$ & $8\ \sC_6$&
$6$\\
\hline
\end{tabular}
\end{center}

Note that only in the last example where $(\alpha,\beta) = (-1,-1)$
is the period of $p_V$ strictly less than $M(\sF)\ D(\sF)$.

We will now justify the entries in the column under $p_V$ in the table.

If $\alpha = \beta =1$, we have solutions
\bigskip
\begin{center}
\begin{tabular}{|l|l|}\hline
$(c,d)$ & solutions $(\theta_1,\theta_2,\theta_3,\theta_4)$\\
\hline
$(1/3,1/3)$&$(1/3,2/3,0,2/3), (5/6,1/6,1/2,2/3)$\\
$(1/3,2/3)$&$(0,2/3,1/3,2/3), (1/2,1/6,5/6,2/3)$\\
$(2/3,1/3)$&$(0,1/3,2/3,1/3), (1/2,5/6,1/6,1/3)$\\
$(2/3,2/3)$&$(2/3,1/3,0,1/3), (1/6,5/6,1/2,1/3)$\\
\hline
\end{tabular}
\end{center}
There are four solutions with order 3 and four with order 6,
thus, 
$$
p_V(n) = \left\{\begin{array}{ll} 8 &\mbox{if $6 \ |\ n$,}\\
4 &\mbox{if $3 \ |\ n$ and $6 \ \not | \ n$}\\
0 &\mbox{otherwise}
\end{array}\right .
$$
and hence $p_V(n) = 4 \ \sC_6(n) + 4 \ \sC_3(n)$.
\smallskip

If $\alpha = 1$ and $\beta = -1$, we have the solutions:
\bigskip
\begin{center}
\begin{tabular}{|l|l|}\hline
$(c,d)$& solutions $(\theta_1,\theta_2,\theta_3,\theta_4)$\\
\hline
$(1/3,1/6)$&$(1/4,11/12,1/12,2/3), (3/4,5/12,7/12,2/3)$\\
$(1/3,5/6)$&$(1/12,5/12,1/4,2/3), (7/12,11/12,3/4,2/3)$\\
$(2/3,1/6)$&$(11/12,7/12,3/4,1/3), (5/12,1/12,1/4,1/3)$\\
$(2/3,5/6)$&$(3/4,1/12,11/12,1/3), (1/4,7/12,5/12,1/3)$\\
\hline
\end{tabular}
\end{center}
All solutions have order $12$ so we have
$p_V(n) = 8 \ \sC_{12}(n)$.
\smallskip

If $\alpha = -1$ and $\beta = 1$, we have the solutions:
\bigskip
\begin{center}
\begin{tabular}{|l|l|}\hline
$(c,d)$& solutions $(\theta_1,\theta_2,\theta_3,\theta_4)$\\
\hline
$(1/6,1/3)$&$(1/4,7/12,11/12,5/6), (3/4,1/12,5/12,5/6)$\\
$(1/6,2/3)$&$(5/12,1/12,3/4,5/6), (11/12,7/12,1/4,5/6)$\\
$(5/6,1/3)$&$(7/12,11/12,1/4,1/6), (1/12,5/12,3/4,1/6)$\\
$(5/6,2/3)$&$(3/4,5/12,1/12,1/6), (1/4,11/12,7/12,1/6)$\\
\hline
\end{tabular}
\end{center}
All solutions have order $12$ so 
$p_V(n) = 8\ \sC_{12}(n)$.

If $a = b = -1$, we have solutions
\bigskip
\begin{center}
\begin{tabular}{|l|l|}\hline
$(c,d)$ & solutions $(\theta_1,\theta_2,\theta_3,\theta_4)$\\
\hline
$(1/6,1/6)$&$(1/6,5/6,0,5/6), (2/3,1/3,1/2,5/6)$\\
$(1/6,5/6)$&$(0,1/3,1/6,5/6), (1/2,5/6,2/3,5/6)$\\
$(5/6,1/6)$&$(0,2/3,5/6,1/6), (1/2,1/6,1/3,1/6)$\\
$(5/6,5/6)$&$(5/6,1/6,0,1/6), (1/3,2/3,1/2,1/6)$\\
\hline
\end{tabular}
\end{center}
All solutions have order 6, 
so $p_V(n) =  8 \ \sC_6(n)$.
\bigskip

\heading{Example 3.}
We now study the Fermat curve $V(f)$, defined by 
$$
f(t_1,t_2) = t_1^m + t_2^m -1.
$$
Then 
\begin{eqnarray*}
&&\Zar{\Tor(V(f))} = \phi_{\sP,f}^{-1}(S(\sP,f))\\
&&S(\sP,f) = \{\mu_1,\mu_2\}
\qquad \mu_1 = (\zeta_{6},\zeta_{6}^{-1}), 
\quad \mu_2 = (\zeta_{6}^{-1},\zeta_{6}),\\
&&M(f) = 6,\\
&&\Epsilon = \Epsilon(\sP,f) = \{(m,0),(0,m)\},\\
&&\overline{\Epsilon} = \{(1,0),(0,1)\} = \Z^2,\\
&&\overline{\Epsilon}/\Epsilon = (\Z/{m\Z})^2,\\
&&D(f) = D((\Z/{m\Z})^2) = m\\
\end{eqnarray*}
We thus have
$$
\deg(p_{V(f)})= 0
$$
and 
$$
6= M(f) \ | \ {\per}(p_{V(f)})\ |\ M(f)\ D(f)= 6m.
$$

For integers $a,b \in \Z$, let $(a,b)$ denote their greatest common divisor.
We'll show that
$$
p_V(f) = \left\{\begin{array}{ll}
2(m,n)^2&\mbox{if $6 \ |\ \frac{n}{(n,m)}$,}\\
0 & \mbox{otherwise.}
\end{array}\right .
$$
and hence $p_{V(f)}$ has period $6m$.  
\smallskip

For instance, if $m = 3$, $p_V(n) = 18\ C_{18}(n)$.
\smallskip

The exponential linear equations associated to
$\phi_{\sP,f}$ are 
\begin{eqnarray}
\left\{\begin{array}{rcl}
m\theta_1 &=& 1/6\ (\mod\ 1)\\
m\theta_2 &=& 5/6\ (\mod\ 1)
\end{array} \right .
\qquad
\left\{\begin{array}{rcl}
m\theta_1 &=& 5/6\ (\mod\ 1)\\
m\theta_2 &=& 1/6\ (\mod\ 1)
\end{array} \right .
\end{eqnarray}

Consider the first equation and suppose there is a solution
$$
(\theta_1, \theta_2) = (a/n,b/n) \in (\Q/\Z)^2. 
$$
We will show that $6$ divides $n/{(m,n)}$.  
The equations in (12) imply that
$$
ma/n - 1/6 \in \Z
$$
so $6n$ divides $6ma - n$ and hence $6$ divides $n$.
Setting $n_1 = n/6$, we have $n$ divides $ma - n_1$.
Set $m_1 = m/(m,n)$.
Then
$$
m_1(m,n)a = n_1\ (\mod\ n).
$$
Since $(m_1,n) = 1$, $m_1$ is invertible modulo $n$ and there is some
$m_2 \in \Z$ such that $m_2m_1 = 1 (\mod\ n)$.
Thus,
$$
(m,n)a = m_2n_1\ (\mod\ n)
$$
and hence
$$
(m,n)a = m_2n_1 + nr,
$$
for some integer $r$.  
Since $n$ is relatively prime
to $m_1$ and $m_2$, so is $(m,n)$.  Thus,
$(m,n)$ divides $n_1$ and hence $6$ divides $n/(m,n)$.

Conversely, suppose $6$ divides $n/{(m,n)}$.  Let $m_1$ be a representative
for the multiplicative inverse of $m/{(m,n)}$ modulo $n$.  Then
$$
a = \frac{n}{6(m,n)}m_1 \qquad b  = \frac{-n}{6(m,n)}m_1
$$
is a solution to the first system and $(b,a)$ is a solution to the second.

To find the number of solutions of order $n$. We need to look 
at the homogeneous equations
$$
m(r/n) = 0\ (\mod\ 1) \qquad m (s/n) = 0\ (\mod\ 1).
$$
These have solutions $r,s\in \Q/\Z$ given by
$$
r = kn/{(m,n)},\quad  s = \ell n/{(m,n)}
$$
for $k,\ell=0,\dots,(m,n) -1$.
Thus, in total, there are $2(m,n)^2$ possible solutions to (12) of order
$n$ when $6$ divides $n/{(m,n)}$.
\qed

\bibliographystyle{math}
\bibliography{math}

\vfill

\hspace{2in}{M.S.R.I., 1000 Centennial Drive, Berkeley, CA 94720}

\hspace{2in}{eko@msri.org}

\vspace{12 pt}
\hspace{2in}{Department of Mathematics, University of Toronto,}

\hspace{2in}{100 St George St., Toronto, ON  M5S 3G3}

\hspace{2in}{eko@math.toronto.edu}

\end{document}